\documentclass[12pt,reqno]{amsart}
\usepackage{amsfonts,amsmath,amssymb}

\allowdisplaybreaks
\advance\textheight by \topskip

\usepackage{mathrsfs}

\usepackage[latin1]{inputenc}
\usepackage{graphicx}

\usepackage{tikz}
\newtheorem{theorem}{Theorem}
\usepackage{euscript}

\def\[{[\! [}
\def\]{]\! ]}

\begin{document}

\title{Skew Dyck paths with catastrophes}

\author[H.~Prodinger]{Helmut Prodinger}

\address{Helmut Prodinger,
		Department of Mathematical Sciences, Stellenbosch University,
7602 Stellenbosch, South Africa, and
 NITheCS (National Institute for Theoretical and Computational Sciences),
South Africa}
\email{hproding@sun.ac.za}

\date{\today}

\begin{abstract}
Skew Dyck paths are like Dyck paths, but an additional south-west step $(-1,-1)$ is allowed, provided that the path does not intersect itself.
Lattice paths with catastrophes can drop from any level to the origin in just one step. We combine these two ideas. The analysis is strictly based
on generating functions, and the kernel method is used.
\end{abstract}

\subjclass{05A15}

\maketitle

\section{Introduction}

The standard random walk on the non-negative integers may be visualized by the following graph (only the first 8 states are shown):

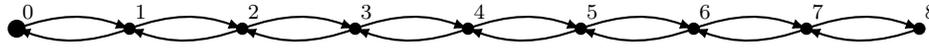
\begin{figure}[h]

	\begin{center}
		\begin{tikzpicture}[scale=1.5]

			\foreach \x in {0,1,2,3,4,5,6,7,8}
			{
				\draw (\x,0) circle (0.05cm);
				\fill (\x,0) circle (0.05cm);
			}

			\fill (0,0) circle (0.08cm);

			\foreach \x in {0,2,4,6}
			{
				\draw[thick,  -latex] (\x,0) to[out=20,in=160]  (\x+1,0);	
				\draw[thick,  -latex] (\x+1,0) to[out=200,in=-20]  (\x,0);	
			}
			\foreach \x in {1,3,5,7}
			{
				\draw[thick,  -latex] (\x,0) to[out=20,in=160]  (\x+1,0);	
				\draw[thick,  -latex] (\x+1,0) to[out=200,in=-20]  (\x,0);	
			}

			\foreach \x in {0,1,2,3,4,5,6,7}
			{
				\node at  (\x+0.1,0.15){\tiny$\x$};
			}

			\node at  (8+0.1,0.15){\tiny$8$};

		\end{tikzpicture}
	\end{center}
	\caption{Standard symmetric random walk on the non-negative integers}
\end{figure}

Such walks are also known as Dyck paths, going up or down one step at the time.
Most of the time, they need to come back to the $x$-axis (end at state 0), but one also
considers open-ended paths where the level of the end of the path is not specified.\footnote{Philippe Flajolet~\cite{BF} liked the names `excursion' resp.\ `meander, which we will not use.}

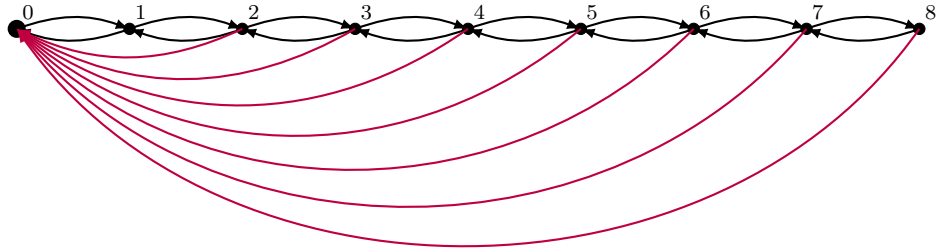
\begin{figure}[h]

	\begin{center}
		\begin{tikzpicture}[scale=1.5]

			\foreach \x in {0,1,2,3,4,5,6,7,8}
			{
				\draw (\x,0) circle (0.05cm);
				\fill (\x,0) circle (0.05cm);
			}

			\fill (0,0) circle (0.08cm);

			\foreach \x in {0,2,4,6}
			{
				\draw[thick,  -latex] (\x,0) to[out=20,in=160]  (\x+1,0);	
				\draw[thick,  -latex] (\x+1,0) to[out=200,in=-20]  (\x,0);	
			}
			\foreach \x in {1,3,5,7}
			{
				\draw[thick,  -latex] (\x,0) to[out=20,in=160]  (\x+1,0);	
				\draw[thick,  -latex] (\x+1,0) to[out=200,in=-20]  (\x,0);	
			}

			\foreach \x in {0,1,2,3,4,5,6,7}
			{
				\node at  (\x+0.1,0.15){\tiny$\x$};
			}

			\node at  (8+0.1,0.15){\tiny$8$};
			
			\draw[thick,purple, -latex] (2,0) to[out=205,in=-25]  (0,0);	
			\draw[thick,purple, -latex] (3,0) to[out=210,in=-30]  (0,0);	
						\draw[thick,purple, -latex] (4,0) to[out=215,in=-35]  (0,0);	
						\draw[thick,purple, -latex] (5,0) to[out=220,in=-40]  (0,0);	
						\draw[thick,purple, -latex] (6,0) to[out=225,in=-45]  (0,0);	
						\draw[thick,purple, -latex] (7,0) to[out=230,in=-50]  (0,0);	
						\draw[thick,purple, -latex] (8,0) to[out=235,in=-55]  (0,0);	
			
		\end{tikzpicture}
		\end{center}
\caption{Catastrophes can lead from a state $i\ge2$ back to the $x$-axis in one step.}
\label{purpel}
\end{figure}
The   Figure \ref{purpel} shows the graph underlying the catastrophes, drawn in purple.

The standard reference for lattice paths with catastrophes is the paper \cite{BW}; the very
recent paper \cite{Baryll} contains some bijective aspects.

The present note will combine skew Dyck paths (with an extra step $(-1,-1)$), as analyzed in
\cite{skew-old}, with the concept of catastrophes.

\section{A warm-up: Dyck paths and catastrophes}

To offer a gentle path for our readers, we analyze the paths as in Fig.~\ref{purpel} (which is also contained in the earlier papers \cite{BW, Baryll}),
since we will approach the skew Dyck paths with catastrophes in a similar style in the next section.

We introduce generating functions $f_i=f_i(z)$, where the coefficient of $z^n$ counts the number of paths starting at the origin (=the big circle) 
and end after $n$ steps at state $i$ (=level $i$). The following recursions are easy to see:
\begin{align*}
f_0&=1+z(f_1+f_2+f_3+f_4+\cdots),\\
f_i&=zf_{i-1}+zf_{i+1}, \ i\ge1.
\end{align*}
Since $f_0$ is somewhat special, we leave it out for the moment and compute the other ones, $f_i$, $i\ge1$.
Eventually we will solve the equation for $f_0$, which will turn out to be just linear.
Therefore we introduce the bivariate generating function
\begin{equation*}
F(u)=F(u,z)=\sum_{i\ge1}u^{i-1}f_{i}
\end{equation*}
and we treat $f_0$ as a parameter. Summing the recursions,
\begin{equation*}
	F(u)=zf_0+zuF(u)+\frac zu[F(u)-f_1],
\end{equation*}
or
\begin{equation*}
F(u)=\frac{zf_0-\frac zuf_1}{1-zu-\frac zu}=\frac{zf_1-zuf_0}{zu^2-u+z}=\frac{zf_1-zuf_0}{z(u-r_1)(u-r_2)},
\end{equation*}
with
\begin{equation*}
r_1=\frac{1+\sqrt{1-4z^2}}{2z},\quad r_2=\frac{1-\sqrt{1-4z^2}}{2z}.
\end{equation*}
An essential step of the \emph{kernel method} is that the `bad' factor $(u-r_2)$ must cancel. The reciprocal of this factor does not
have a power series expansion in $(u,z)$. Dividing this factor out in both, denominator and numerator, leads to
\begin{equation*}
	F(u)=\frac{-zf_0}{z(u-r_1)}=\frac{zf_0}{zr_1(1-u/r_1)},
\end{equation*}
and from this
\begin{equation*}
f_i=[u^{i-1}]F(u)=\frac{zf_0}{zr_1^i}=f_0r_2^i.
\end{equation*}
This formula holds for all $i\ge0$. It is interesting to note that all these functions are just a multiple of the parameter $f_0$. Now we can go to the first recursion:
\begin{equation*}
f_0=1+z(f_0r_2^1+f_0r_2^2+f_0r_2^3+\cdots)=1+zf_0\frac{r_2}{1-r_2}
\end{equation*}
and solve the equation,
\begin{equation*}
f_0=\frac{2-3z-2z^2+z\sqrt{1-4z^2}}{2(1-z-2z^2-z^3)}.
\end{equation*}
It has the power series expansion 
\begin{equation*}
f_0=1+z^2+z^3+3z^4+5z^5+12z^6+23z^7+52z^8+105z^9+232z^{10}+480z^{11}+1049z^{12}+2199z^{13}+\cdots,
\end{equation*}
and the list of coefficients is sequence A224747 in the OEIS \cite{OEIS}. The papers \cite{BW, Baryll} have this already. We can also compute the number of Dyck paths 
with catastrophes and open (unspecified) end, via the generating function
\begin{align*}
f_0&+f_1+f_2+f_3+f_4+\dots\\
&=\frac{1-z+(1+z)\sqrt{1-4z^2}}{2(1-z-2z^2-z^3)}\\
&=1+z+2{z}^{2}+4{z}^{3}+8{z}^{4}+17{z}^{5}+35{z}^{6}+75{z}^
{7}+157{z}^{8}+337{z}^{9}+712{z}^{10}+\cdots.
\end{align*}
This is almost sequence A274115 in the OEIS \cite{OEIS}; $1+z(f_0+f_1+f_2+\cdots)$ has this sequence as coefficients. Thanks go to Michel Marcus to point out the 
connection with the sequences in the OEIS.

\section{Skew Dyck paths with catastrophes}
Skew Dyck are a variation of Dyck paths, where additionally to steps $(1,1)$ and $(1,-1)$ a south-west step $(-1,-1)$ is also allowed, provided that the path
does not intersect itself.

Here is a list of the 10 skew paths consisting of 6 steps:

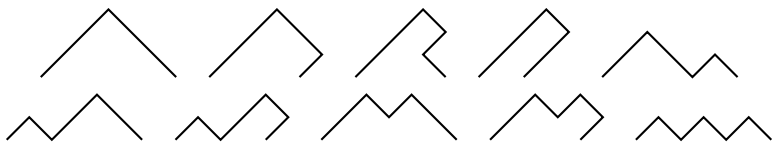
\begin{figure}[h]
	\begin{equation*}
		\begin{tikzpicture}[scale=0.3]
			\draw [thick](0,0)--(3,3)--(6,0);
		\end{tikzpicture}
		\quad
		\begin{tikzpicture}[scale=0.3]
			\draw [thick](0,0)--(3,3)--(5,1)--(4,0);
		\end{tikzpicture}
		\quad
		\begin{tikzpicture}[scale=0.3]
			\draw [thick](0,0)--(3,3)--(4,2)--(3,1)--(4,0);
		\end{tikzpicture}
		\quad
		\begin{tikzpicture}[scale=0.3]
			\draw [thick](0,0)--(3,3)--(4,2)--(3,1)--(2,0);
		\end{tikzpicture}
		\quad
		\begin{tikzpicture}[scale=0.3]
			\draw [thick](0,0)--(2,2)--(4,0)--(5,1)--(6,0);
		\end{tikzpicture}
	\end{equation*}
	\begin{equation*}
		\begin{tikzpicture}[scale=0.3]
			\draw [thick](0,0)--(1,1)--(2,0)--(4,2)--(6,0);
		\end{tikzpicture}
		\quad
		\begin{tikzpicture}[scale=0.3]
			\draw [thick](0,0)--(1,1)--(2,0)--(4,2)--(5,1)--(4,0);
		\end{tikzpicture}
		\quad
		\begin{tikzpicture}[scale=0.3]
			\draw [thick](0,0)--(1,1)--(2,2)--(3,1)--(4,2)--(6,0);
		\end{tikzpicture}
		\quad
		\begin{tikzpicture}[scale=0.3]
			\draw [thick](0,0)--(1,1)--(2,2)--(3,1)--(4,2)--(5,1)--(4,0);
		\end{tikzpicture}
		\quad
		\begin{tikzpicture}[scale=0.3]
			\draw [thick](0,0)--(1,1)--(2,0)--(3,1)--(4,0)--(5,1)--(6,0);
		\end{tikzpicture}
	\end{equation*}
	\caption{All 10 skew Dyck paths of length 6 (consisting of 6 steps).}
\label{list1}
\end{figure}

We prefer to work with the equivalent model (resembling more traditional Dyck paths) where
we replace each step $(-1,-1)$ by $(1,-1)$ but label it red. Here is the list of the 10 paths again (Figure \ref{list2}):

\begin{figure}[h]
	\begin{equation*}
		\begin{tikzpicture}[scale=0.3]
			\draw [thick](0,0)--(3,3)--(6,0);
		\end{tikzpicture}
		\quad
		\begin{tikzpicture}[scale=0.3]
			\draw [thick](0,0)--(3,3)--(5,1);
			\draw [thick,red](5,1)--(6,0);
		\end{tikzpicture}
		\quad
		\begin{tikzpicture}[scale=0.3]
			\draw [thick](0,0)--(3,3)--(4,2);
			\draw[red,thick] (4,2)--(5,1);
			\draw [thick](5,1)--(6,0);
		\end{tikzpicture}
		\quad
		\begin{tikzpicture}[scale=0.3]
			\draw [thick](0,0)--(3,3)--(4,2);
			\draw[red,thick](4,2)--(6,0);
		\end{tikzpicture}
		\quad
		\begin{tikzpicture}[scale=0.3]
			\draw [thick](0,0)--(2,2)--(4,0)--(5,1)--(6,0);
		\end{tikzpicture}
	\end{equation*}
	\begin{equation*}
		\begin{tikzpicture}[scale=0.3]
			\draw [thick](0,0)--(1,1)--(2,0)--(4,2)--(6,0);
		\end{tikzpicture}
		\quad
		\begin{tikzpicture}[scale=0.3]
			\draw [thick](0,0)--(1,1)--(2,0)--(4,2)--(5,1);
			\draw[red,thick] (5,1)--(6,0);
		\end{tikzpicture}
		\quad
		\begin{tikzpicture}[scale=0.3]
			\draw [thick](0,0)--(1,1)--(2,2)--(3,1)--(4,2)--(6,0);
		\end{tikzpicture}
		\quad
		\begin{tikzpicture}[scale=0.3]
			\draw [thick](0,0)--(1,1)--(2,2)--(3,1)--(4,2)--(5,1);
			\draw[red,thick] (5,1)--(6,0);
		\end{tikzpicture}
		\quad
		\begin{tikzpicture}[scale=0.3]
			\draw [thick](0,0)--(1,1)--(2,0)--(3,1)--(4,0)--(5,1)--(6,0);
		\end{tikzpicture}
	\end{equation*}
	\caption{The 10 paths redrawn, with red south-east edges instead of south-west edges.}
\label{list2}
\end{figure}
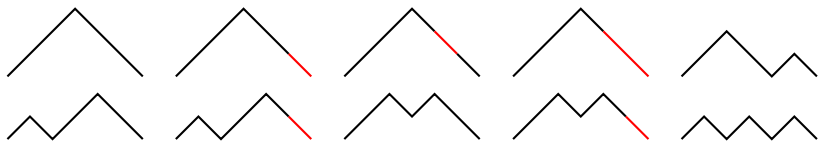

The rules to generate such decorated Dyck paths are: each edge $(1,-1)$ may be  black or red, but
\begin{tikzpicture}[scale=0.3]\draw [thick](0,0)--(1,1); \draw [red,thick] (1,1)--(2,0);\end{tikzpicture}
and
\begin{tikzpicture}[scale=0.3] \draw [red,thick] (0,1)--(1,0);\draw [thick](1,0)--(2,1);\end{tikzpicture}
are forbidden.

These figures are taken from our earlier paper \cite{skew-old}.

\begin{figure}[h]

	\begin{center}
		\begin{tikzpicture}[scale=1.5]
			\draw (0,0) circle (0.1cm);
			\fill (0,0) circle (0.1cm);
			
			\foreach \x in {0,1,2,3,4,5,6,7,8}
			{
				\draw (\x,0) circle (0.05cm);
				\fill (\x,0) circle (0.05cm);
			}
			
			\foreach \x in {0,1,2,3,4,5,6,7,8}
			{
				\draw (\x,-1) circle (0.05cm);
				\fill (\x,-1) circle (0.05cm);
			}
			
			\foreach \x in {0,1,2,3,4,5,6,7,8}
			{
				\draw (\x,-2) circle (0.05cm);
				\fill (\x,-2) circle (0.05cm);
			}
			
			\foreach \x in {0,1,2,3,4,5,6,7}
			{
				\draw[ thick,-latex] (\x,0) -- (\x+1,0);
				
			}

			\foreach \x in {1,2,3,4,5,6,7}
			{
				\draw[thick,  -latex] (\x+1,0) to[out=200,in=70]  (\x,-1);

			}
			\draw[ thick,  -latex] (1,0) to[out=200,in=70]  (0,-1);

			\foreach \x in {0,1,2,3,4,5,6,7}
			{
				
				\draw[thick,  -latex] (\x,-1) to[out=30,in=250]  (\x+1,0);	
				
			}

			\foreach \x in {0,1,2,3,4,5,6,7}
			{
				\draw[ thick,-latex] (\x+1,-1) -- (\x,-1);
				
			}
			\foreach \x in {0,1,2,3,4,5,6,7}
			{
				\draw[ thick,-latex,red] (\x+1,-1) -- (\x,-2);
				
			}
			
			\foreach \x in {0,1,2,3,4,5,6,7}
			{
				\draw[ thick,-latex,red] (\x+1,-2) -- (\x,-2);
				
			}
			
			\foreach \x in {0,1,2,3,4,5,6,7}
			{
				\draw[ thick,-latex] (\x+1,-2) -- (\x,-1);
				
			}
			
			\foreach \x in {2,3,4,5,6,7,8}
			{
				\draw[ ultra thick, purple,-latex] (\x,0) -- (\x-0.25,0.25);
				
			}
			\foreach \x in {2,3,4,5,6,7,8}
		{
			\draw[ ultra thick, purple,-latex] (\x,0-1) -- (\x-0.25,0.25-1);
			
		}
			\foreach \x in {2,3,4,5,6,7,8}
		{
			\draw[ ultra thick, purple,-latex] (\x,0-2) -- (\x-0.25,-0.25-2);
			
		}
	
		\end{tikzpicture}
	\end{center}
	\caption{Three layers of states according to the type of steps leading to them (up, down-black, down-red). Purple arrows all end in the origin.}
	\label{threelayers}
\end{figure}
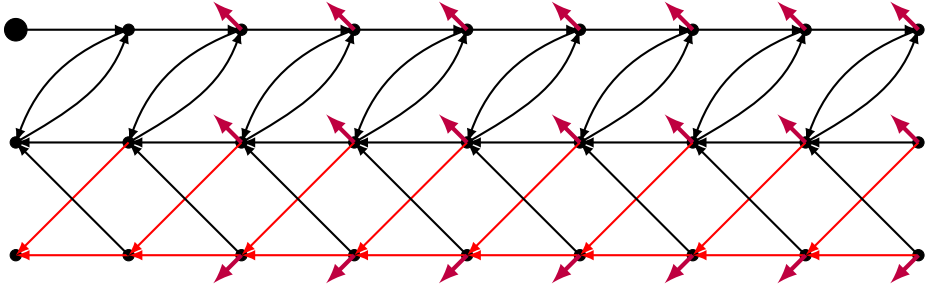

As in the motivating example from the previous section, we introduce generating functions. But since we need 3 layers to control everything, 
we need to introduce $f_j$, $g_j$, $h_j$. The purple edges all end in the origin and represent the catastrophes. If they are ignored, one models
skew Dyck paths in this way, in particular $f_0+g_0+h_0$ is the generating function of skew Dyck paths coming back to level 0.
However, now, we will deal with the purple arrows as well.

The following recursions can be read off immediately from the diagram \ref{threelayers}:
\begin{gather*}
	f_0=1+z\sum_{i\ge2}f_i+z\sum_{i\ge2}g_i+z\sum_{i\ge2}h_i,\\
	f_{i+1}=zf_i+zg_i,\quad i\ge0,\\
	g_i=zf_{i+1}+zg_{i+1}+zh_{i+1},\quad i\ge0,\\
	h_i=zh_{i+1}+zg_{i+1},\quad i\ge0.
\end{gather*}
As before, we first ignore $f_0$ and treat is as a parameter. We introduce
\begin{equation*}
F(u)=\sum_{i\ge1}u^{i-1}f_i,\quad G(u)=\sum_{i\ge0}u^{i}g_i,\quad H(u)=\sum_{i\ge0}u^{i}h_i.
\end{equation*}
Now we sum the recursions and get
\begin{align*}
	\sum_{i\ge0}u^if_{i+1}=\sum_{i\ge0}u^izf_i+\sum_{i\ge0}u^izg_i
\end{align*}
or
\begin{align*}
	F(u)=f_0+uzF(u)+zG(u).
\end{align*}
Further,
\begin{equation*}
	\sum_{i\ge0}u^ig_i=\sum_{i\ge0}u^izf_{i+1}+\sum_{i\ge0}u^izg_{i+1}+\sum_{i\ge0}u^izh_{i+1}
\end{equation*}
which translates into
\begin{equation*}
	G(u)=zF(u)+\frac zu[G(u)-G(0)]+\frac zu[H(u)-H(0)];
\end{equation*}
similarly
\begin{equation*}
	H(u)=\frac zu[G(u)-G(0)]+\frac zu[H(u)-H(0)].
\end{equation*}
We will eliminate the functions $G$ and $H$: 
\begin{equation*}
G=-{\frac {z{ f_0}+zuF-F}{z}}
\end{equation*}
and
\begin{equation*}
	H=G-zF.
\end{equation*}
Therefore we end up with just one equation for $F$
\begin{equation*}
-uzf_0-z{u}^{2}F+uF=-{z}^{2}Fu-2 {z}^{2}f_0+2 zF-2g_0{z}^{2}-{z}^{3}F+{z}^{3}f_1
\end{equation*}
which can be solved,
\begin{equation*}
F(u)=-{\frac { \left( -2zf_0-2 zg_0+{z}^{2}f_1+f_0
		u \right) z}{-{z}^{2}u+2z-{z}^{3}+z{u}^{2}-u}}
=-\frac { \left( -2zf_0-2 zg_0+{z}^{2}f_1+f_0
		u \right) z}{z(u-r_1)(u-r_2)}
\end{equation*}
with
\begin{equation*}
	r_1=\frac{1+z^2+\sqrt{1-6z^2+5z^4}}{2z},\quad r_2=\frac{1+z^2-\sqrt{1-6z^2+5z^4}}{2z}.
\end{equation*}
The fact that $u-r_2$ must divide the numerator of $F$ leads to
\begin{equation*}
g_0={\frac {-2 zf_0+{z}^{2}f_1+r_2f_0}{2z}},
\end{equation*}
so that $g_0$ can be eliminated and only $f_1$ remains as unknown;
\begin{equation*}
F=\frac{zf_0}{1-zu+z^2-zr_2}
\end{equation*}
plugging in $u=0$ allows to compute it:
\begin{equation*}
f_1=\frac{zf_0}{1+z^2-zr_2}=\frac{f_0}{r_1}.
\end{equation*}
Therefore the function $F$ is known:
\begin{equation*}
F(u)=\frac{f_0}{r_1-u}
\end{equation*}
and thus
\begin{equation*}
f_j=[u^{j-1}]F(u)=\frac{f_0}{r_1^j}.
\end{equation*}
From this, the other functions can be computed as well:
\begin{equation*}
G(u)=\frac{1-zr_1}{z}\frac{f_0}{r_1-u}=(r_2-z)\frac{f_0}{r_1-u}
\end{equation*}
and 
\begin{equation*}
g_j=[u^{j}]G(u)=(r_2-z)\frac{f_0}{r_1^{j+1}},
\end{equation*}
as well as
\begin{equation*}
H(u)=(r_2-2z)\frac{f_0}{r_1-u}
\end{equation*}
and
\begin{equation*}
	h_j=[u^{j}]H(u)=(r_2-2z)\frac{f_0}{r_1^{j+1}}.
\end{equation*}
Now we can finally deal with $f_0$. The evaluations that are needed are simple:
\begin{equation*}
\sum_{j\ge2}f_j=\frac{f_0/r_1^2}{1-1/r_1},\quad
\sum_{j\ge2}g_j=(r_2-z)\frac{f_0/r_1^3}{1-1/r_1},\quad
\sum_{j\ge2}h_j=(r_2-2z)\frac{f_0/r_1^3}{1-1/r_1};
\end{equation*}
of course this is better done with a computer. 
From this, the equation (just linear) for $f_0$ can be solved:
\begin{equation*}
f_0=\frac{1-9 {z}^{3}-5 {z}^{4}+3 {z}^{5}+2 {z}^{6}-(1-z^2-z^3)\sqrt{1-6z^2+5z^4}}{2(1-2 {z}^{2}-6 {z}^{3}-3 {z}^{4}+ {z}^{5}+ {z}^{6})}.
\end{equation*}
The series expansion
\begin{equation*}
1+{z}^{3}+{z}^{4}+4 {z}^{5}+6 {z}^{6}+18 {z}^{7}+31 {z}^{8}+85 {
	z}^{9}+157 {z}^{10}+410 {z}^{11}+792 {z}^{12}+2004 {z}^{13}+\cdots
\end{equation*}
is not in OEIS \cite{OEIS}. To compute the number of all skew Dyck paths coming back to the $x$-axis, we need
\begin{align*}
f_0&+g_0+h_0\\*&=
\frac{(1-z^2)(3-z)(1-z-2z^2-z^3)-(1-4z-3z^2+z^3+z^4)\sqrt{1-6z^2+5z^4}}{2(1-2 {z}^{2}-6 {z}^{3}-3 {z}^{4}+ {z}^{5}+ {z}^{6})}.
\end{align*}
The series expansion is
\begin{equation*}
	1+z^2+z^3+4z^4+5z^5+17z^6+25z^7+76z^8+125z^9+353z^{10}+625z^{11}+1681z^{12}+\cdots.
\end{equation*}

One can also sum everything to enumerate the open ended paths in this model:
\begin{align*}
\sum_{j\ge0}&f_j+\sum_{j\ge0}g_j+\sum_{j\ge0}h_j\\
&=\frac{(1+z)(1-4z+4z^2+4z^3-z^5)+(1+5z+3z^2-z^3-z^4)\sqrt{1-6z^2+5z^4}}{2(1-2{z}^{2}-6{z}^{3}-3{z}^{4}+{z}^{5}+{z}^{6})};
\end{align*}
the series expansion is
\begin{equation*}
1+z+2z^2+4z^3+9z^4+18z^5+41z^6+85z^7+193z^8+410z^9+929z^{10}+2004z^{11}+\cdots.
\end{equation*}
\begin{theorem}
	The generating functions related to the model of skew Dyck paths with catastrophes are
	\begin{align*}
f_j&=\frac{f_0}{r_1^j},\\
g_j&=(r_2-z)\frac{f_0}{r_1^{j+1}},\\
h_j&=(r_2-2z)\frac{f_0}{r_1^{j+1}},
	\end{align*}
with 
\begin{equation*}
r_1=\frac{1+z^2+\sqrt{1-6z^2+5z^4}}{2z},\quad r_2=\frac{1+z^2-\sqrt{1-6z^2+5z^4}}{2z},
\end{equation*}
and
\begin{align*}
	f_0&+g_0+h_0\\*&=
	\frac{(1-z^2)(3-z)(1-z-2z^2-z^3)-(1-4z-3z^2+z^3+z^4)\sqrt{1-6z^2+5z^4}}{2(1-2 {z}^{2}-6 {z}^{3}-3 {z}^{4}+ {z}^{5}+ {z}^{6})}.
\end{align*}
The enumeration of all paths, regardless of the final level, is
\begin{align*}
	\sum_{j\ge0}&f_j+\sum_{j\ge0}g_j+\sum_{j\ge0}h_j\\
	&=\frac{(1+z)(1-4z+4z^2+4z^3-z^5)+(1+5z+3z^2-z^3-z^4)\sqrt{1-6z^2+5z^4}}{2(1-2{z}^{2}-6{z}^{3}-3{z}^{4}+{z}^{5}+{z}^{6})}.
\end{align*}

\end{theorem}

\section{Conclusion}

From our analysis, other types of catastrophes (example: only from an even-numbered level one can drop to the ground) can be computed without much further effort. This changes only the equation for $f_0$; the other quantities stay as they are.


\end{document}